\documentclass[11pt]{article}
\usepackage{amscd,amssymb,epic,amsmath}

\usepackage[dvips]{graphicx}
\usepackage[mathscr]{eucal}
\usepackage[all]{xy}
\usepackage{eepic}

\newtheorem{lemma}{Lemma}

\newtheorem{theorem}{Theorem}
\newtheorem{remark}{Remark}

\title{\Large \bf On $*$-representations of polynomial algebras in quantum matrix spaces of rank 2}

\author{Olga Bershtein\\
Tallinn University of Technology,\\
e-mail:olga.bershteyn@ttu.ee}

\date{}

\begin{document}

\maketitle

\begin{abstract}
In this paper we study of $*$-representations for polynomial algebras on quantum matrix spaces. We deal with two special cases of the polynomial algebras, namely the algebra of polynomials on quantum complex matrices $\mathrm{Mat_2}$ and on quantum complex symmetric matrices $\mathrm{Mat_2^{sym}}$. For the second algebra we classify all irreducible $*$-representations by bounded operators in a Hilbert space (up to a unitary equivalence). Moreover, we present a construction of $*$-representations of the above algebras which enables to obtain the full list of $*$-representations (sometimes by passing to subrepresentations).
\end{abstract}

\medskip


\medskip
\section{Introduction}

The theory of quantum groups established many rich topics during the last 30 years. One of them is the theory of quantum bounded symmetric domains. These domains were introduced and investigated by L.~Vaksman and his collaborators (see \cite{SV} and the monograph \cite{V} together with the references therein). Recall that Harish-Chandra suggested the standard realization of a classical bounded symmetric domain $\mathbb D$ as a unit ball in a certain normed vector space $\mathfrak{p}$. Thus the main object in the quantum theory is a $*$-algebra of polynomials $\mathrm{Pol}(\mathfrak{p})_q$ endowed with a full $U_q \mathfrak{g}$-symmetry, where $\mathfrak{g}=\mathfrak{g_0}\otimes_{\mathbb R} \mathbb C, \mathfrak{g_0}=Lie(G_0)$ and $G_0$ is the automorphism group of $\mathbb D$. These algebras lead to geometric realizations for representations of quantum groups, noncommutative version for some problems of complex analysis, quantum Harish-Chandra modules, harmonic analysis related to these domains and special functions.

D.~Shklyarov and L.~Vaksman suggested a way to construct $*$-representations (by bounded operators) of polynomial algebras in complex matrix spaces. This approach enables to obtain a list of $*$-representations from certain simplest representations (namely, faithfull representations of such algebras) and well known representations of $\mathbb C[U_n]_q$ (see \cite{Soib}). It is worth to say that the construction is nearly automatic. But then one has to check the irreducibility for such representations and the completeness of the list.

The last problem naturally belongs to a wide topic of studying representations for $*$-algebras and its classification. There is a vast literature on this subject. Let me mention just a survey \cite{OS}. In this book $*$-representations are classified by means of studying the corresponding dynamical systems.

In this paper we deal with matrix spaces that are assigned to bounded symmetric domains of rank 2. Under Harish-Chandra realization, these domains are just unit balls in the spaces of complex and symmetric complex $2 \times 2$-matrices. Thus the polynomial algebras in quantum bounded symmetric domains can be considered as $*$-algebras in the respective quantum vector spaces $\mathrm{Mat_2}$ and $\mathrm{Mat_2^{sym}}$. In the first case $*$-representations of $\mathrm{Pol}(Mat_2)_q$ by bounded operators in a Hilbert space were classified by L.Turowska in \cite{T}. In Appendix we just present our construction of $*$-representations for this special case. For $\mathrm{Pol}(Mat_2^{\rm sym})_q$ we obtain a list of $*$-representations and verify its completeness. Let us also mention that our proof is similar to the proof of the main theorem in \cite{T}.

In what follows $\mathbb C$ is the ground field, all algebras are associative and unital, $q \in (0,1)$.
\bigskip

\section{$*$-Representations of $\mathrm{Pol}(Mat_2^{\rm sym})_q$}

Let us recall the required notion and notation in our special case. The bounded symmetric domain of type $C_2$ admits a standard realization as a unit ball $\mathbb D$ in the space of symmetric $2 \times 2$-matrices. In the quantum case all necessary definitions and notation were introduced in \cite{B1} for a general bounded symmetric domain of type $C_n$. Here we just specify them to our needs.

Now the $*$-algebra $\mathrm{Pol(Mat_2^{\rm sym})}_q$ can be defined by generators $z_{11},z_{21},z_{22}$ and the following list of relations:
\begin{align*}
& z_{11}z_{21}=q^2z_{21}z_{11}, \quad z_{21}z_{22}=q^2z_{22}z_{21},\\
& z_{11}z_{22}-z_{22}z_{11}=q(q^2-q^{-2})z_{21}^2,\\
& z_{11}^*z_{11}=q^4z_{11}z_{11}^*-q(q^{-1}-q)(1+q^2)^2z_{21}z_{21}^*+
(q^{-1}-q)^2(1+q^2)z_{22}z^*_{22}+1-q^4,\\
& z_{11}^*z_{21}=q^2z_{21}z_{11}^*-q(q^{-1}-q)(q^{-1}+q)z_{22}z_{21}^*\\
& z_{11}^*z_{22}=z_{22}z^*_{11}, \quad z_{21}^*z_{22}=q^2z_{22}z_{21}^*,\\
& z_{21}^*z_{21}=q^2z_{21}z_{21}^*-(1-q^2)z_{22}z_{22}^*+(1-q^2),\\
& z_{22}^*z_{22}=q^4z_{22}z_{22}^*+1-q^4.
\end{align*}

\begin{remark}
To make this definition more symmetric, one may introduce an additional generator $z_{12}$ together with the relation $z_{12}=qz_{21}$.
\end{remark}

Let us obtain a list of all irreducible $*$-representations of $\mathrm{Pol(Mat_2^{sym})}_q$ by bounded operators in a Hilbert space. For that we use the standard technique of analysing the spectrum of a commutative subalgebra and the related dynamical system \cite{OS}.

It follows from the commutation relations that $z_{21}z_{21}^*, z_{22}z_{22}^*$ form a commutative subalgebra in $\mathrm{Pol(Mat_2^{sym})}_q$. Moreover, one can compute that for $(a,b)=(2,1)$ or $(a,b)=(2,2)$ and each $(c,d)$
\begin{equation*}
z_{ab}z_{ab}^*z_{cd}=z_{cd}F_{b}^{d}(z_{21}z_{21}^*,z_{22}z_{22}^*),
\end{equation*}
where
\begin{align*}
& F_{1}(x_1,x_2)=(F_1^1,F_1^2)=(q^2x_1-(1-q^2)x_2+1-q^2,x_2),\\
& F_{2}(x_1,x_2)=(F_2^1,F_2^2)=(q^4x_1,q^4x_2+1-q^4).
\end{align*}
Let us denote by $F$ the action of $\mathbb Z^2$ on $\mathbb R^2$ induced by these functions. The corresponding orbits are
\begin{equation*}
\Omega_{x_1,x_2}=\{(q^{2m}q^{4n}x_1-q^{4n}(1-q^{2m})(x_2-1),q^{4n}x_2+1-q^{4n})|m,n \in \mathbb Z\}.
\end{equation*}

Let $\pi$ be an irreducible $*$-representation of $\mathrm{Pol(Mat_2^{sym})}_q$ in a Hilbert space $H$ by bounded operators. Denote by $E$ the resolution of identity for the commutative system $\pi(z_{21}z_{21}^*), \pi(z_{22}z_{22}^*)$. One can prove similarly to Lemma 1 of \cite{T} that

\begin{lemma}
There exist $x_1,x_2 \in \mathbb R_+$ such that $E(\Omega_{x_1,x_2})=I$.
\end{lemma}
{\bf Proof.} Let us sketch the proof. If a subset $\Delta \subset \mathbb R^2$ is invariant under $F$ then it induces a subspace $E(\Delta) \subset H$, which is invariant under $\pi_{z_{21}},\pi_{z_{21}^*},\pi_{z_{22}},\pi_{z_{22}^*}$. Moreover, since
\begin{align*}
& (z_{21}z_{21}^*)z_{11}=z_{11}z_{21}z_{21}^*+q(q^2-q^{-2})z_{21}^2z_{22}^*,\\
& (z_{22}z_{22}^*)z_{11}=z_{11}z_{22}z_{22}^*-q(q^2-q^{-2})z_{21}^2z_{22}^*,
\end{align*}
$E(\Delta)$ is invariant under $\pi(z_{11}),\pi(z_{11}^*)$ either. Thus $E$ is ergodic with respect to the dynamical system generated by the $F$-action. This action is one-to-one and there exists a Borel set which meets every orbit $\Omega_{x_1,x_2}$ in a single point. Thus we have the statement of the lemma. \hfill $\square$

Now let us prove the following
\begin{theorem}
Each irreducible $*$-representation of $\mathrm{Pol(Mat_2^{sym})}_q$ by bounded operators in a Hilbert space is unitary equivalent to a representation from the following list
\begin{enumerate}
\item A two-parameter series of representations in $\mathbb C$
\begin{align*}
& \pi^{(1)}_{\phi,\psi}(z_{11})=e^{i\psi},
\\ & \pi^{(1)}_{\phi,\psi}(z_{21})=0, \quad \pi^{(1)}_{\phi,\psi}(z_{22})=e^{i\phi}, \qquad \psi,\phi \in \mathbb R/2\pi \mathbb Z,
\end{align*}
\item A one-parameter series of representations in $l_2(\mathbb Z_+)$
\begin{align*} & \pi^{(2)}_{\phi}(z_{11})e_k=q^{-1}\sqrt{1-q^{4k+4}}e_{k+1},
\\ & \pi^{(2)}_{\phi}(z_{21})e_k=0, \quad \pi^{(2)}_{\phi}(z_{22})e_k=e^{i\phi}e_k, \qquad k \in \mathbb Z_+, \phi \in \mathbb R/2\pi \mathbb Z,
\end{align*}
\item A one-parameter series of representations in $l_2(\mathbb Z_+)$
\begin{align*}
& \pi^{(3)}_{\phi}(z_{11})e_k=-q^{-1}\sqrt{1-q^{4k}}e_{k-1},
\\& \pi^{(3)}_{\phi}(z_{21})e_k=q^{2k}e^{i\phi}e_k,
\\& \pi^{(3)}_{\phi}(z_{22})e_k=\sqrt{1-q^{4k+4}}e_{k+1}, \qquad k \in \mathbb Z_+, \phi \in \mathbb R/2\pi \mathbb Z.
\end{align*}
\item A one-parameter series of representations in $l_2(\mathbb Z_+)^{\otimes 2}$
\begin{align*}
& \pi^{(4)}_{\phi}(z_{11})e_{k,l}=q^{2k}e^{i\phi}e_{k,l}-q^{-1}\sqrt{(1-q^{4l})(1-q^{2k+2})(1-q^{2k+4})}e_{k+2,l-1},
\\& \pi^{(4)}_{\phi}(z_{21})e_{k,l}=q^{2l}\sqrt{1-q^{2k+2}}e_{k+1,l},
\\& \pi^{(4)}_{\phi}(z_{22})e_{k,l}= \sqrt{1-q^{4l+4}}e_{k,l+1}, \qquad k,l \in \mathbb Z_+, \phi \in \mathbb R/2\pi \mathbb Z,
\end{align*}
\item A representation in $l_2(\mathbb Z_+)^{\otimes 3}$
\begin{align*}
& \pi^{(5)}(z_{11})e_{k,l,m}=q^{2l}\sqrt{1-q^{4m+4}}e_{k,l,m+1}-q^{-1}\sqrt{(1-q^{4k})(1-q^{2l+2})(1-q^{2l+4})}e_{k-1,l+2,m},
\\& \pi^{(5)}(z_{21})e_{k,l,m}=q^{2k}\sqrt{1-q^{2l+2}}e_{k,l+1,m},
\\& \pi^{(5)}(z_{22})e_{k,l,m}= \sqrt{1-q^{4k+4}}e_{k+1,l,m}, \qquad k,l,m \in \mathbb Z_+.
\end{align*}
\end{enumerate}
\end{theorem}

{\bf Proof.} It remains to analyze which orbits $\Omega_{x_1,x_2}$ lead to irreducible representations by bounded operators in a Hilbert space. An argument quite similar to those from \cite{T} proves that only $\Omega_{0,0}, \Omega_{1,0}, \Omega_{0,1}$ give bounded operators. Let us consider these cases one by one.

Denote by $\sigma$ the joint spectrum of $\pi(z_{21}z_{21}^*), \pi(z_{22}z_{22}^*)$. Let $\sigma \subset \Omega_{0,1}$. One can easily show that $\pi(z_{21})=0,\pi(z_{22})=e^{i\phi}{\rm id}$. Now for the remaining generator $z_{11}$ we have
\begin{equation*}
\pi(z_{11}^*)\pi(z_{11})=q^4\pi(z_{11})\pi(z_{11}^*)+q^{-2}-q^2.
\end{equation*}
This relation (after renormalization) leads to the standard one-dimensional $q$-oscillator algebra. For this algebra irreducible $*$-representations by bounded operators in a Hilbert space are well known (see, for example, \cite[Section 1.4]{OS}). The list of representations consists of an infinite-dimensional representation in $l_2(\mathbb Z_+)$ and a series of one-dimensional representations. Thus we obtain the series $\pi^{(1)}_{\phi,\psi}$ and $\pi^{(2)}_{\phi}$.

Let $\sigma \subset \Omega_{1,0}=\{(q^{4n}, 1-q^{4n})|n \in \mathbb Z\}$. Put
\begin{equation*}
H_n=\{v \in H|\pi(z_{21})\pi(z_{21}^*)v=q^{4n}v,\pi(z_{22})\pi(z_{22}^*)v=(1-q^{4n})v\}.
\end{equation*}

Using the commutation relations, one can verify that $\pi(z_{21})(H_n) \subset H_n$ and $\pi(z_{22})(H_n) \subset H_{n+1}$. Denote by $P_n$ the projection to $H_n$ parallel to other eigenspaces. Since
\begin{equation*}
(z_{21}z_{21}^*)z_{11}=z_{11}(z_{21}z_{21}^*)+q(q^-2-q^{-2})z_{21}^2z_{22}^*,
\end{equation*}
one has that
\begin{equation*}
\pi(z_{11})P_n=P_n\pi(z_{11})P_n+P_m\pi(z_{11})P_n, \qquad m=n-1.
\end{equation*}
Moreover,
\begin{equation*}
P_{n-1}\pi(z_{11})P_n=q\frac{q^2-q^{-2}}{q^{4n-4}(1-q^4)}z_{21}^2z_{22}^*P_n=
\\ -q^3P_{n-1}z_{21}^2z_{22}^*(1-\pi(z_{22}z_{22}^*))^{-1}P_n.
\end{equation*}
Thus the action of $\pi(z_{11})$ splits into the sum of its diagonal part $\pi(z_{11})_0$ and the operator $-q^3z_{21}^2z_{22}^*(1-\pi(z_{22}z_{22}^*))^{-1}$.
By explicit computations,
\begin{equation*}
\pi(z_{11})^*_0\pi(z_{11})_0=q^4\pi(z_{11})_0\pi(z_{11})_0^*.
\end{equation*}
Therefore, $\pi(z_{11})_0$ has to be zero operator. So the operator of representation $\pi(z_{11})$ is written by means of $\pi(z_{21})$ and $\pi(z_{22})$, so we may not take it into account in our irreducibility considerations. This yields that each $H_n$ has to be irreducible under the action of $\pi(z_{21})$. Moreover, if we consider the polar decomposition $\pi(z_{21})=UP$, the same will hold for the unitary part $U$. On the other hand, $U$ commutes with $\pi(z_{22}), \pi(z_{22})^*$, thus it is just a scalar operator in $H$ and each $H_n$ is one-dimensional. Now we may consider the linear span of the vectors $\pi(z_{22})^kv,v \in H_0$, and this gives us $\pi^{(3)}_{\phi}$.

Now let us consider the case $\sigma \subset \Omega_{0,0}=\{(q^{4n}(1-q^{2m}), 1-q^{4n})|m,n \in \mathbb Z\}$. Similarly, put
\begin{equation*}
H_{m,n}=\{v \in H|\pi(z_{21})\pi(z_{21}^*)v=q^{4n}(1-q^{2m})v,\pi(z_{22})\pi(z_{22}^*)v=(1-q^{4n})v\}.
\end{equation*}
Using the commutation relations, one can verify that $\pi(z_{21})(H_{m,n}) \subset H_{m+1,n}$ and $\pi(z_{22})(H_{m,n}) \subset H_{m,n+1}$. Besides, $\pi(z_{11})(H_{m,n}) \subset H_{m,n}\oplus H_{m+2,n-1}$. In this case the diagonal part of $\pi(z_{11})$ quasicommutes with all other generators and satisfies the following relation (here $P_{m,n}$ denotes the projection to $H_{m,n}$ parallel to other eigenspaces):
\begin{equation*}
\pi(z_{11})_0^*\pi(z_{11})_0P_{m,n}=q^4\pi(z_{11})_0\pi(z_{11})_0^*P_{m,n}+(1-q^4)q^{4m}.
\end{equation*}
Thus $\pi(z_{11})_0$ satisfies the standard relation in $H_{0,0}$ and the action of $\pi(z_{11})_0,\pi(z_{11})_0^*$ has to be simple here. It leads to either the infinite dimensional representation of the $q$-oscillator algebra in $H_{0,0}$ or one of its one dimensional representations. On this way we obtain $\pi^{(4)}_{\phi}$ and $\pi^{(5)}$ $*$-representations. \hfill $\square$

\bigskip

Now we present a construction of $*$-representations of $\mathrm{Pol(Mat_2^{sym})}_q$.

It is well-known that $\mathbb D$ is a homogeneous space of the group $U_2$ under the action
\begin{equation}\label{act_sym}
U: Z \mapsto UZU^t, \qquad U \in U_2, Z \in \mathbb D \subset \mathrm{Mat_2^{sym}}.
\end{equation}
Let us restrict ourselves by the $SU_2$-action for a while. Consider the following $*$-homomorphism of coaction
\begin{align*}
\Delta: \mathrm{Pol(Mat_2^{sym})}_q \rightarrow \mathrm{Pol(Mat_2^{sym})}_q
\otimes \mathbb C[SU_2]_q,
\\ \Delta(z_{jk})= z_{11} \otimes t_{1j} t_{1k}+
qz_{21} \otimes t_{1j} t_{2k} + z_{21} \otimes t_{2j} t_{1k}+ z_{22} \otimes
t_{2j} t_{2k}.
\end{align*}

This homomorphism is extensively used in our construction of $*$-representations of $\mathrm{Pol(Mat_2^{sym})}_q$.

Now let us introduce the simplest representations of $\mathrm{Pol(Mat_2^{sym})}_q$.

It is proved in \cite{V} that for any quantum bounded symmetric domain $\mathbb D$ there exists a unique (up to a unitary equivalence) irreducible $*$-representation of $\mathrm{Pol}(\mathfrak{p})_q$ by bounded operators in a Hilbert space. The corresponding module is defined by a single generator $v$ and the relations
\begin{equation*}
z_{11}^*v=z_{21}^*v=z_{22}^*v=0.
\end{equation*}
This representation is called the Fock representation. Let us denote by $F_2$ the Fock representation of $\mathrm{Pol(Mat_2^{sym})}_q$.

It follows from the commutation relations that one has the $*$-morhpism of algebras
$$
\pi_1 :\mathrm{Pol(Mat_2^{sym})}_q \rightarrow \mathrm{Pol(\mathbb
C)}_{q^2}, \quad \pi_1
\begin{pmatrix}
z_{11} & qz_{21}\\ z_{21} & z_{22}
\end{pmatrix}
\mapsto \begin{pmatrix} q^{-1} z & 0 \\ 0 & 1
\end{pmatrix}.
$$
Here the $*$-algebra $\mathrm{Pol(\mathbb
C)}_{t}$ is defined by a single generator $z$ and the well-known relation $z^*z=t^2zz^*+1-t^2$. This algebra naturally corresponds to the simplest bounded symmetric domain, namely the unit disc $\{z \in \mathbb C| |z|<1\}$.

Put $F_1=T \circ \pi_1$, where $T$ is the Fock representation of $\mathrm{Pol(\mathbb C )}_{q^2}$ (again, the corresponding module is defined by a generator $v$ and the relation $z^*v=0$).

Now we consider the $*$-morphism
$$
F_0 :\mathrm{Pol(Mat_2^{sym})}_q \rightarrow \mathbb C, \quad F_0
\begin{pmatrix}
z_{11} & qz_{21}\\ z_{21} & z_{22}
\end{pmatrix}
\mapsto \begin{pmatrix} q^{-1}  & 0 \\ 0 & 1
\end{pmatrix}.
$$

It is well-known (see \cite{Soib}) that $\mathbb C[SU_2]_q$ possesses an infinite dimensional $*$-representation $\pi: \mathbb C[SU_2]_q \rightarrow \mathrm{End} (l_2(\mathbb Z_+))$ which in the standard basis $\{e_n\},n \in \mathbb Z_+$ is defined as follows:
\begin{align}\label{repSU2_1}
\pi(t_{11})e_n= \sqrt{1-q^{2n}}e_{n-1}, & \qquad \pi(t_{12})e_n=q^{n+1}e_n,
\\ \pi(t_{21})e_n=-q^n e_n, & \qquad \pi(t_{22})e_n=\sqrt{1-q^{2n+2}}e_{n+1}, \nonumber
\end{align}
and the one-dimensional $*$-representation
\begin{equation}\label{repSU2_2}
\varepsilon:\mathbb C[SU_2]_q \rightarrow \mathbb C, \quad
\varepsilon(t_{ij})=\delta_{ij}.
\end{equation}

Now let us turn to a construction of $*$-representations. Let $\rho_1$ be one of the simplest representations of $\mathrm{Pol(Mat_2^{sym})}_q$ and $\rho_2$ one of the simplest representations of $\mathbb C[SU_2]_q$ described above. We consider representations of $\mathrm{Pol(Mat_2^{sym})}_q$ of the form $(\rho_1\boxtimes \rho_2)\Delta$ (we use the notation $\boxtimes$ since our factors are representations of different algebras). Thus we get the following list:
\begin{enumerate}
\item $(F_2 \boxtimes \pi)\Delta$. One can show that this representation is reducible. Namely, it has a subrepresentation isomorphic to $F_2$.

\item $(F_2 \boxtimes \varepsilon)\Delta$. It is the Fock representation itself.

\item $(F_1 \boxtimes \pi)\Delta$. This representation has a natural realization in $l_2(\mathbb Z_+)^{\otimes 2}$. Indeed, if we put $v=e_0 \boxtimes e_0 \in l_2(\mathbb Z_+)^{\otimes 2}$ then the action of generators is defined as follows:
$$\begin{array}{lll}z_{11} v=v,&& \\
(z_{11})^*v=v,& (z_{21})^*v=(z_{22})^*v=0.\end{array}$$
One can easily check that this representation is irreducible.

\item $(F_1 \boxtimes \varepsilon)\Delta$. This representation has a natural realization in $l_2(\mathbb Z_+)$. If we put $v=e_0 \in l_2(\mathbb Z_+)$ then the action of generators is defined as follows:
$$\begin{array}{lllll}& z_{21} v=0,& z_{22} v=v,& \\
(z_{11})^*v=0,& (z_{21})^*v=0,& (z_{22})^*v=v.\end{array}$$
This representation is also irreducible.

\item $(F_0 \boxtimes \pi)\Delta$. This representation admits a realization in $l_2(\mathbb Z_+)$.
An easy calculation shows that this representation is reducible. Namely, it contains two subrepresentations isomorphic to $\pi^{(3)}_0$.

\item $(F_0 \boxtimes \varepsilon)\Delta$. This one-dimensional representation is defined via
$$\begin{array}{lllll}z_{11} =q^{-1},& z_{21} v=0,& z_{22} =1,&
\\ (z_{11})^*=q^{-1},& (z_{21})^*=0,& (z_{22})^*=1.\end{array}$$
\end{enumerate}

Now we can recall an additional parameter. Indeed, one has an action of the torus $\mathbb T^2$ in $Mat_2^{sym}$ given also by \eqref{act_sym}. This action can be naturally transferred to representations of $\mathrm{Pol(Mat_2^{sym})}_q$. Thus we obtain the following list of $*$-representations of $\mathrm{Pol(Mat_2^{sym})}_q$:
\begin{enumerate}
\item The vector space of the representation $(F_2 \boxtimes \varepsilon)\Delta$ is spanned by $z_{22}^az_{21}^bz_{11}^cv, a,b,c \in \mathbb Z_+$. One can introduce a scalar product and complete the vector space with respect to it. Thus we get $\pi^{(5)}$ from Theorem 1.

\item The one-parameter series of $*$-representations described as follows:
$$\begin{array}{lll}z_{11} v=e^{i\psi}v,&& \\
(z_{11})^*v=e^{-i\psi}v,& (z_{21})^*v=(z_{22})^*v=0, & \psi \in
\mathbb{R}/2 \pi \mathbb{Z}.\end{array}$$
The vector space of this representation is spanned by $z_{21}^az_{22}^bv, a,b \in \mathbb Z_+$. After completion with respect to suitable scalar product we get $\pi^{(4)}_{\psi}$ from Theorem 1.

\item The one-parameter series of $*$-representations described as follows:
$$\begin{array}{lllll}& z_{21} v=0,& z_{22} v=e^{i\psi}v,& \\
(z_{11})^*v=0,& (z_{21})^*v=0,& (z_{22})^*v=e^{-i\psi}v, & \psi \in
\mathbb{R}/2 \pi \mathbb{Z}.\end{array}$$
It is exactly the $\pi^{(1)}_{\phi}$ series of representations (after a completion with respect to suitable scalar product).

\item The one-parameter series of $*$-representations described as follows:
$$\begin{array}{llll}
& z_{21}v=e^{i\phi}v,
\\ (z_{11})^*v=q^{-1}e^{2i\phi}z_{22} v,& (z_{21})^*v=-e^{-i\phi}v,
\,\, (z_{22})^*v=0, & \phi, \in \mathbb{R}/2 \pi
\mathbb{Z}.\end{array}$$
This series can be obtained from $(F_0 \boxtimes \pi)\Delta$ by passing to a subrepresentations. After a suitable completion, one would obtain $\pi^{(3)}_{\phi}$.

\item The two-parameter series $\pi^{(2)}_{\phi,\psi}$ of $*$-representations described as follows:
$$\begin{array}{lllll}z_{11} =q^{-1}e^{i\phi},& z_{21} v=0,&
z_{22} =e^{i\psi},& \\ (z_{11})^*=q^{-1}e^{-i\phi},& (z_{21})^*=0,&
(z_{22})^*=e^{-i\psi}, & \varphi, \psi \in \mathbb{R}/2 \pi
\mathbb{Z}.\end{array}$$
\end{enumerate}

\section{Appendix. A construction of $*$-representations of $\mathrm{Pol}(Mat_2)_q$}

In this appendix we describe a list of irreducible $*$-representations of $\mathrm{Pol(Mat_2)}_q$ using our constructions. This list was previously obtained by Vaksman, Sinel'shchikov and Shklyarov in \cite{SSV2} and by Turowska in \cite {T}, where its completeness was proved.

Recall that $\mathrm{Pol(Mat_2)}_q$ can be defined in terms of generators $z_1^1,z_1^2,z_2^1,z_2^2$ and the following relations:
\begin{align*}
& z_i^jz_k^l=qz_k^lz_i^j \qquad \qquad i=k \,\& \,j \neq l \quad \text{or} \quad i \neq k \,\& \,j=l,
\\ & z_1^2z_2^1=z_2^1z_1^2,
\\ & z_1^1z_2^2=z_2^2z_1^1+(q-q^{-1})z_1^2z_2^1,
\end{align*}
\begin{align*}
& (z_1^1)^*z_1^1=q^2z_1^1(z_1^1)^*-q^2(q^{-2}-1)(z_2^1(z_2^1)^*+z_1^2(z_1^2)^*)+q^2(q^{-2}-1)^2z_2^2(z_2^2)^*+1-q^2,
\\ & (z_1^1)^*z_i^j=qz_i^j(z_1^1)^*+(q-q^{-1})z_2^2(z_i^j)^*, \qquad \qquad (i,j)=(1,2) \quad \text{or} \quad (i,j)=(2,1),
\\ & (z_1^1)^*z_2^2=z_2^2(z_1^1)^*,
\\ & (z_i^j)^*z_i^j=q^2z_i^j(z_i^j)^*-(1-q^2)z_2^2(z_2^2)^*+1-q^2, \qquad \qquad (i,j)=(1,2) \quad \text{or} \quad (i,j)=(2,1),
\\ & (z_2^1)^*z_1^2=z_1^2(z_2^1)^*,
\\ & (z_2^2)^*z_i^j=qz_i^j(z_2^2)^*, \qquad \qquad (i,j)=(1,2) \quad \text{or} \quad (i,j)=(2,1),
\\ & (z_2^2)^*z_2^2= q^2z_2^2(z_2^2)^*+1-q^2.
\end{align*}

Let us consider the $*$-homomorphism of coaction
$$ \mathcal{D}: \mathrm{Pol(Mat_2)}_q \rightarrow \mathrm{Pol(Mat_2)}_q \otimes
\mathbb C[SU_2 \times SU_2]_q, \qquad \mathcal{D}(z_j^i)=\sum z_b^a \otimes t_{bj}
\otimes t_{ai}.
$$
Recall that this coaction corresponds to well known left-right action of $SU_2 \times SU_2$ in $\mathrm{Mat_2}$
\begin{equation}\label{act_mat}
(A,D): Z \mapsto AZD^{-1}, \qquad A,D \in SU_2, Z \in \mathrm{Mat_2}.
\end{equation}
Our construction involves only the representations obtained through $\mathcal{D}$.

Let us introduce the simplest known $*$-representations of $\mathrm{Pol(Mat_2)}_q$.

Firstly, one has the Fock representation $\mathcal{F}_2$ of $\mathrm{Pol(Mat_2)}_q$ (i.e. a faithfull irreducible $*$-representation by bounded operators in a Hilbert space, see \cite{SSV4}). The corresponding module is defined by a single generator $v$ and the relations
\begin{equation*}
(z_i^j)^*v=0, \qquad i,j=1,2.
\end{equation*}

By obvious reasons, one has the $*$-morphism of algebras
$$
\pi_1 :\mathrm{Pol(Mat_2)}_q \rightarrow \mathrm{Pol(\mathbb C)}_q, \quad
\pi_1
\begin{pmatrix}
z_1^1 & z_1^2\\ z_2^1 & z_2^2
\end{pmatrix}
\mapsto \begin{pmatrix} q^{-1} z & 0 \\ 0 & 1
\end{pmatrix}.
$$

Now let $\mathcal{F}_1=T \circ \pi_1$ be the second simplest representation of $\mathrm{Pol(Mat_2)}_q$ (here $T$ is the Fock representation of $\mathrm{Pol(\mathbb C )}_q$).

Thirdly, define the $*$-homomorphism
$$
\mathcal{F}_0 :\mathrm{Pol(Mat_2)}_q \rightarrow \mathbb C, \quad \mathcal{F}_0
\begin{pmatrix}
z_1^1 & z_1^2\\ z_2^1 & z_2^2
\end{pmatrix}
\mapsto \begin{pmatrix} q^{-1}  & 0 \\ 0 & 1
\end{pmatrix}.
$$

Now, using the standard representations of $\mathbb C[SU_2]_q$ \eqref{repSU2_1}-\eqref{repSU2_2}, we get the following list:
\begin{enumerate}
\item $(\mathcal{F}_2 \boxtimes \pi \boxtimes\pi)\mathcal{D}$.

\item $(\mathcal{F}_2 \boxtimes \pi \boxtimes \varepsilon)\mathcal{D}$ or $(\mathcal{F}_2 \boxtimes \varepsilon
\boxtimes \pi)\mathcal{D}$.

\item $(\mathcal{F}_2 \boxtimes \varepsilon \boxtimes \varepsilon)\mathcal{D}$. It is isomorphic to the Fock representation of $\mathrm{Pol(Mat_2)}_q$ itself.

\item $(\mathcal{F}_1 \boxtimes \pi \boxtimes \pi)\mathcal{D}$. It has a natural realization $l_2(\mathbb Z_+)^{\otimes 3}$. Put $v=e_0 \boxtimes e_0 \boxtimes e_0$ then one has
$$\begin{array}{lll}z_1^1 v=v,&& \\
(z_1^1)^*v=v,& (z_1^2)^*v=(z_2^1)^*v=(z_2^2)^*v=0.\end{array}$$

\item $(\mathcal{F}_1 \boxtimes \pi \boxtimes \varepsilon)\mathcal{D}$ or $(\mathcal{F}_1 \boxtimes \varepsilon
\boxtimes \pi)\mathcal{D}$. These representations can be realized in $l_2(\mathbb Z_+)^{\otimes 2}$ as follows:
$v= e_0 \boxtimes e_0$ and
$$\begin{array}{lllll} z_1^1 v=0,& & z_2^1 v=v,&& \\
(z_1^1)^*v=0,& (z_1^2)^*v=0,& (z_2^1)^*v=v,& (z_2^2)^*v=0.\end{array}
$$
or
$$\begin{array}{lllll} z_1^1 v=0,& & z_2^1 v=v,&& \\
(z_1^1)^*v=0,& (z_1^2)^*v=0,& (z_2^1)^*v=v,& (z_2^2)^*v=0.\end{array}$$

\item $(\mathcal{F}_1 \boxtimes \varepsilon \boxtimes \varepsilon)\mathcal{D}$. It admits a natural realization in $l_2(\mathbb Z_+)$: put $v=e_0$, then one has
$$\begin{array}{lllll}& z_1^2 v=0,& z_2^1 v=0,& z_2^2 v=v,& \\
(z_1^1)^*v=0,& (z_1^2)^* v=0,& (z_2^1)^*v=0,& (z_2^2)^*v=v.\end{array}$$

\item $(\mathcal{F}_0 \boxtimes \pi \boxtimes \pi)\mathcal{D}$.

\item $(\mathcal{F}_0 \boxtimes \pi \boxtimes \varepsilon)\mathcal{D}$ or $(\mathcal{F}_0 \boxtimes \varepsilon
\boxtimes \pi)\mathcal{D}$. These representations are isomorphic and have natural realizations in $l_2(\mathbb Z_+)$: put $v=e_0$ then one has
$$\begin{array}{llll}z_1^1 v=0,& z_1^2 v=v,& z_2^1 v=v,&
\\ (z_1^1)^*v=-q^{-1}z_2^2 v,& (z_1^2)^* v=v,&
(z_2^1)^*v=v,& (z_2^2)^*v=0.\end{array}$$

\item $(\mathcal{F}_0 \boxtimes \varepsilon \boxtimes \varepsilon\mathcal{D}$. This one-dimensional representation is defined by the following:
$$\begin{array}{lllll}z_1^1
=q^{-1},& z_1^2 =0,& z_2^1 v=0,& z_2^2 =1,& \\ (z_1^1)^*=q^{-1},& (z_1^2)^*
=0,& (z_2^1)^*=0,& (z_2^2)^*=1.\end{array}$$
\end{enumerate}

One can verify that the cases 1,2, and 7 give reducible representations. Indeed, in the first two cases the representations have a subrepresentation equivalent to the Fock representation. An explicit calculations in the seventh case shows that this representation have subrepresentations equivalent to $\rho_{\phi,\pi-\phi}$, $\phi \in [0,2\pi)$ (see the defining relations below). Namely, one can directly check that $(z_2^2)^*$ has a countably dimensional kernel, $(z_2^1)^*=-z_1^2$ and they both preserve $\rm{Ker} (z_2^2)^*$. Moreover, for any $\phi \in [0,2\pi)$ there exists $v \in Ker (z_2^2)^*$ such that $z_2^1v=e^{i\phi}v$.

In all other cases the representations are irreducible.

Similar to the case of quantum symmetric matrices, we can recall the action of the torus $\mathbb T^4$ in $\mathbb D \subset \mathrm{Mat_2}$ which corresponds to \eqref{act_mat}. Thus we obtain the following list of irreducible $*$-representations of $\mathrm{Pol(Mat_2)}_q$ (here I use the notation of L.~Turowska):
\begin{enumerate}
\item A faithfull representation $\rho$ of $\mathrm{Pol(Mat_2)}_q$ in $l_2(\mathbb Z_+)^{\otimes 4}$. This representation can be obtained from $(\mathcal{F}_2 \boxtimes \varepsilon \boxtimes \varepsilon)\mathcal{D}$ by completion with respect to a suitable scalar product.

\item A one-parameter series of $*$-representations $\check{\rho}_\phi$, $\phi \in [0,2\pi)$ in $l_2(\mathbb Z_+)^{\otimes 3}$. These representations can be obtained from $(\mathcal{F}_1 \boxtimes \pi \boxtimes \pi)\mathcal{D}$ by completion.

\item Two one-parameter series of $*$-representations $\rho^{(1)}_\phi$ and $\rho^{(2)}_\phi$, $\phi \in [0,2\pi)$ in $l_2(\mathbb Z_+)^{\otimes 2}$. These representations can be obtained from $(\mathcal{F}_2 \boxtimes \pi \boxtimes \varepsilon)\mathcal{D}$ or $(\mathcal{F}_2 \boxtimes \varepsilon \boxtimes \pi)\mathcal{D}$ by completion.

\item A one-parameter series of $*$-representations $\pi_\phi$, $\phi \in [0,2\pi)$ in $l_2(\mathbb Z_+)$. These representations can be obtained from $(\mathcal{F}_1 \boxtimes \varepsilon \boxtimes \varepsilon)\mathcal{D}$ by completion.

\item A two-parameter series of $*$-representations $\rho_{\phi_1,\phi_2}$, $\phi_1,\phi_2 \in [0,2\pi)$ in $l_2(\mathbb Z_+)$.
    These representations can be obtained from $(\mathcal{F}_0 \boxtimes \pi \boxtimes \varepsilon)\mathcal{D}$ by completion.

\item A two-parameter series of one dimensional $*$-representations $\xi_{\phi_1,\phi_2}$, $\phi_1,\phi_2 \in [0,2\pi)$.
\end{enumerate}

\end{document}